\documentclass[journal]{IEEEtran}
\usepackage{amsmath, amsthm, amssymb,graphicx,bbm}
\usepackage{verbatim}
\usepackage{psfrag}
\usepackage{algorithm}
\usepackage{algorithmic}
\usepackage{bm}
\usepackage{color}
\usepackage[tight,footnotesize]{subfigure}

\newcommand{\bd}{\mathbf}

\newcommand{\bma}{\begin{bmatrix}}
\newcommand{\ebma}{\end{bmatrix}}

\newcommand{\mc}{\mathcal}

\DeclareMathOperator{\E}{\mathbb{E}}

\DeclareMathOperator{\cov}{cov}

\newtheorem{defn}{Definition}

\usepackage{amsmath, amsthm, amssymb,graphicx,bbm}

\title{Competition and Coalition Formation of \\Renewable Power Producers}
\author{Baosen Zhang, Ramesh Johari, Ram Rajagopal
\thanks{All authors are at Stanford University. B. Zhang is with the Departments of Civil and Environmental Engineering and Management Sciences and Engineering. R. Johari is with the Department of Management Sciences and Engineering. R. Rajagopal is with the Department of Civil and Environmental Engineering.  Emails: \{zhangbao,ramesh.johari,ramr\}@stanford.edu 

 This work is partially supported by the Stanford Precourt Institute for Energy Efficiency and the US-Israel Binational Foundation.}}

\begin{document}

\maketitle

\begin{abstract}
We investigate group formations and strategic behaviors of renewable power producers in electricity markets. These producers currently bid into the day-ahead market in a conservative fashion because of the real-time risk associated with not meeting their bid amount. It has been suggested in the literature that producers would bid less conservatively if they can form large groups to take advantages of spatial diversity to reduce the uncertainty in their aggregate output. We show that large groups of renewable producers would act strategically to lower the aggregate output because of market power. To maximize renewable power production, we characterize the trade-off between market power and generation uncertainty as a function of the size of the groups. We show there is a sweet spot in the sense that there exists groups that are large enough to achieve the uncertainty reduction of the grand coalition, but are small enough such that they have no significant market power. We consider both independent and correlated forecast errors under a fixed real-time penalty. We also consider a real-time market where both selling and buying of energy are allowed.  We validate our claims using PJM and NREL data.
\end{abstract}
\section{Introduction}
Renewable resources such as wind and solar are expected to play increasingly prominent roles in power systems. An aspect of the bulk electricity system that is fundamental to the success of integration of these renewable sources is the electricity market. Since it is the main venue of resource allocation in power systems \cite{Kirschen04}, understanding the interaction between market rules and the producers is crucial to maximize the gain from renewables. Because the key difference between renewable resources and traditional generators is that the former is much more \emph{uncertain} than the latter,  the main question of interest is how does \emph{uncertainty} impact the outcome and function of the market. 


Currently, most electricity markets in the United States and Europe operate in a multistage manner \cite{Kirschen04,Pinson14}. Twenty four hours in advance of the actual operating time (often termed real-time),  a day-ahead market is used to match supply and demand for an hour long slot. Then one or more additional stages are used to adjust for variations in supply and demand that may not have been settled in the day-ahead market. However, the current market was not designed to accommodate large amounts of uncertainty in the supply.  A guiding philosophy of electricity market design is to ensure that supply and demand are always balanced, so that generators are incentivized to meet their day-ahead allocations. Generators that cannot meet their obligations are subjected to a penalty on their shortfall, or extremely volatile real-time prices, or a combination of both \cite{Pinson14}. 


Because deviations between real-time output and day-ahead obligations are disincentivized, renewable producers tend to bid conservatively in electricity markets to protect against real-time risks. Since the day-ahead prediction error of a wind farm can be up to 25\% (somewhat less for solar), producers bid much less than the forecasted amount of renewables \cite{GE10,NERC09,Holttinen11}. This behavior in turn limits the actual amount of power generated from renewables, since the load not served by the renewables is picked up by traditional generators \cite{Varaiya11}. 



In addition to improving the forecasting technology, one promising method of reducing the uncertainty in renewable resources is to take advantage of geographical diversity as pointed out by many authors \cite{NERC09,EnerNex10,Diakov12,Zhao13,GQ12,Makarov09,Pinson13}. For example, aggregating renewable producers at spatially separate locations can reduce the variability of the total output. In essence, the aggregate is easier to forecast than its individual parts, thus an aggregate of producers could bid less conservatively into the day-ahead market and consequently increase the amount of power generated from renewable sources \cite{Klessmann08,Abbad10,Morales10,Bitar12,Baeyens13}. Therefore it seems that system operators should encourage renewable producers to form coalitions.  

\begin{figure}[ht]
	\centering
	\includegraphics[scale=0.5]{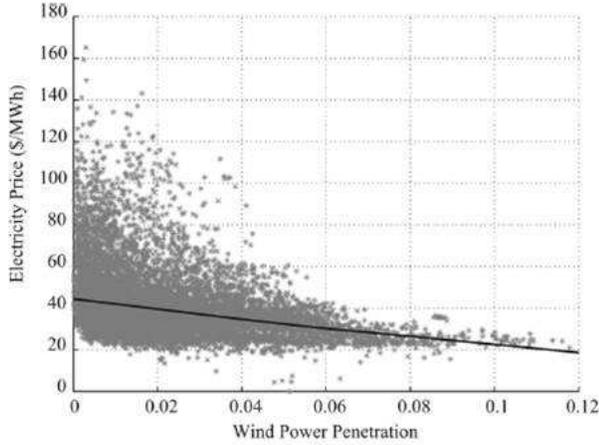}
	\caption{Scatter plot of PJM day-ahead prices and wind generation for 2012. Figure reproduced from \cite{Gil13}. The horizontal axis is the percentage penetration of wind, and the vertical axis is the average clearing price in the PJM area. As the amount of wind penetration increases from 0 to 10\%, the average price in the system drops by more than half.} 
	\label{fig:pjm}
	\end{figure}

On the other hand, aggregating or grouping renewable producers could potentially increase their market power and lead to unintended outcomes. Because of several high profile cases\footnote{For example, see the California Electricity Crisis \cite{FERC}.}, system operators vigilantly oppose generator collusion of any kind. Indeed,  if an aggregate of renewable producers is large enough, it could become price-anticipatory since its bid could have significant effects on the clearing prices of the market. Fig \ref{fig:pjm} (reproduced from  \cite{Gil13}) shows the day-ahead prices drop dramatically even at 10\% penetration of wind in the PJM control area\footnote{PJM is a regional transmission organization covering parts of the eastern United States.}. In order to maximize its profit (quantity times price), an aggregate may not bid all of its forecasted wind. Therefore an aggregate of renewable producers may withhold some of the renewable not because of uncertainty, but because of the market power. Thus from this point of view, it seems that operators should not allow producers to form groups. 
	
In this paper, we show that the system operator can have its cake and eat it too. More precisely, we investigate the trade-off between uncertainty and market power for an aggregation of renewable power producers. We show that for a wide range of scenarios, there exists groups of certain sizes that \emph{both} induce the maximum amount of renewable penetration and do not possess significant market power. Interestingly, the grand coalition (group of all producers) is never desirable.  To arrive at this conclusion, we first build a parametrized model of the residual demand curve to isolate the effect of renewables on the day-ahead clearing price. Each group of renewable producers is allowed to bid a number that represent their production quantity into the market. The resulting product of the clearing price and the bid quantity determines the day-ahead profit of the groups of producers. For the real-time risk, we consider two models: a fixed penalty on the shortfall and a stochastic model for the real-time market. The payoff for each group is defined as the sum of the day-ahead profit plus the real-time risk, and we set up a Cournot game based on the payoff functions. We establish the trade-off between uncertainty and market power by  analyzing the Nash equilibrium of this game. Some of the more game-theoretic questions are studied in our paper in \cite{Zhang14}. 

In this paper (and \cite{Zhang14}) the uncertainty in the renewable production can be thought as follows. Each producer has an estimate of its available power (e.g. from wind or solar irradiation) in the form of a random variable. Note the randomness comes from the error associated with the estimate.  For a group, the estimate is the sum of the estimates of its members, which is again a random variable. First we consider the case when the random variables are independent (that is, the forecast error of the producers are independent). We scale the system by letting the number of producers go to infinity, and the size of each producer go to zero, so the total amount of renewable generation in the system is fixed. By the law of large numbers, as long as the number of producers in a group grows, the uncertainty of the aggregate goes to zero (random variable for that group concentrates around its mean). Also, it is widely known that as the number of players grows, Cournot games becomes competitive since no player has significant market power \cite{Johari05}. Our result essentially states that as long as the number of groups and the size of each group both go to infinity, the market is competitive (competition between groups) and the uncertainties are mitigated (averaging within a group). We extend this idea to the case where producers' estimates are correlated.  We note that this paper focuses on the empirical behaviors of the renewable producers and the corresponding consequences for the electricity market. For a more in depth study of the Cournot game, interested readers can refer to \cite{Zhang14}.

Coalitions of wind producers have been considered in the past by many authors \cite{Zhao13,Pinson13,Bitar12,Baeyens13,Nayyar13}. However the main focus was on how to divide up the profit of a group among its members, not on the effect of groups on the entire system. Since most of previous studies assume that wind farms are always price takers, the grand coalition is often the most desirable set up. Recently, some authors have focused on the strategic behavior of wind farms \cite{Nieta14,Sharma14}, which is closer to our setting. But these papers mainly consider the strategic action of a single price-anticipatory wind farm, whereas we focus on the joint behavior of many producers. 

This paper is organized as follows. Section II introduces the model and problem setup. Section III studies the effect of coalition under a fixed real-time penalty for both independent and correlated producers. Section IV investigates the effect of a real-time market on the risk faced by the producers.  Finally, Section V concludes the paper.

\section{Model and Problem Setup}
We consider a system consisting of renewable producers, traditional generators, and loads. Throughout, $N$ represents the number of renewable producers. We do not distinguish between traditional generators. We assume the loads in the system are inelastic, deterministic and known to all parties. We will assume that there is no congestion and use a single bus model for the network. The case of a network with possible congestions is not considered in this paper and it is an important direction of future research. We use the terms group and coalition interchangeably in this paper. 

We adopt a two-stage structure for the delivery of electricity consisting of a day-ahead stage followed by a real-time stage. The day-ahead stage is a pool-based market, where generators submit their bids to a system operator. The operator clears the auction and determines the generation schedule. The clearing price is denoted by $p$ (\$/MW). If a generator cannot meet its promised amount, a penalty is assessed for the shortfall at real-time.  In some markets, generators are assessed the real-time cost (or profit) for the net deviation at the real-time profit. We consider the fixed penalty case in Section \ref{sec:penalty} and real-time market case in Section \ref{sec:real-time}.  It is well know that the real-time prices are notoriously hard to predict and model \cite{Astaneh13,Allcott11}, and a simple stochastic model is adopted in Section \ref{sec:real-time}.


Since renewable producers have zero marginal cost (or near zero cost), we restrict them to bid only the \emph{amount} of energy they are willing to deliver into the day-ahead market\footnote{Traditional generators typically bid a curve representing the cost of generation at different output levels.}.  Let $w_i$ be the bid by producer and $W_i$ be the random variable representing the amount of renewable generation. Note $W_i$ is based on the day-ahead forecast information, so the randomness in $W_i$ can be thought of as coming from the \emph{forecast error}. Under the constant real-time penalty,  the total expected payoff for producer $i$ is denoted $\pi_i$, given by
\begin{equation}\label{eqn:pi}
\pi_i(w_i)=p(w_1,\dots,w_N)w_i - q\E[(w_i-W_i)^+]. 
\end{equation} 
where $(\cdot)^+$ is the positive part of a number. The first term reflects that the bids from all producers affect the day-ahead price $p$.  The second term reflect the penalty term: if the realized renewable power, $W_i$, is less than the amount bid, $w_i$, a penalty is paid based on the real-time price. Here, the coefficient $q$ can be thought as a given constant\footnote{Equivalently, $q$ can be thought as a random variable that is independent to everything else, the we may replace $q$ by $\E[q]$.}. In Section \ref{sec:real-time}, we consider a real-time market in which the second term can be both positive or negative.

\subsection{Impact of Renewables on the Day-ahead Price}
In this section we determine the impact of renewables on day-ahead market clearing prices.  Since the demand is inelastic and we assume that the traditional generation mix is fixed, the clearing price without any renewable injection can be normalized to be $1$ per unit. As renewable producers bid into the market, the clearing price would drop below that. Also, we normalize the demand to be $1$ per unit since it is inelastic.  

\begin{figure}[ht]
\hspace{-0.3cm}
\includegraphics[scale=0.5]{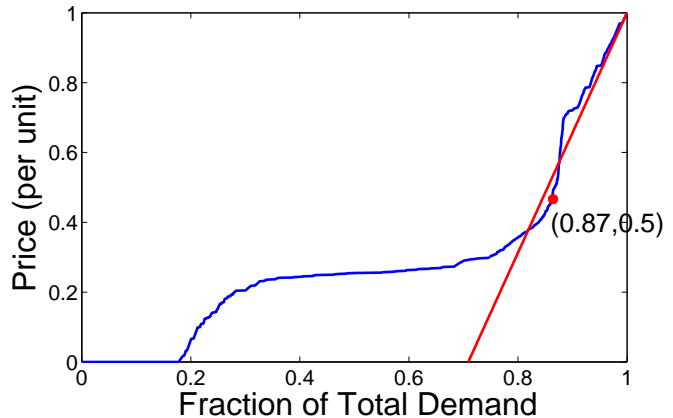}
\caption{ The demand curve for PJM 2007 bids. The red dot shows that about 87\% of the total demand can be purchased at half of the clearing price. For the data used, the clearing price was at \$60.
The red line shows the linear approximation for the right part of the demand curve. This approximation works well up to about 25\% renewable penetration.}
\label{fig:pjm_bid}
\end{figure}

The demand curve is defined as a function relating the clearing price of purchasing one unit of energy to the total amount of energy being purchased. With the above normalizations, the residual demand curve is defined on the interval $[0,1]$  and takes values on the interval $[0,1]$. To construct the demand curve for a particular day, the bid curves of all generators would be stacked up to determine the cheapest clearing price for a certain demand. Since we are interested on the overall qualitative behavior of renewable producers, we study the average demand curve constructed from historical bid information. Fig. \ref{fig:pjm_bid} shows the demand curve for PJM in 2007. This year was chosen to capture the generator bids before there was significant wind in the system. 


The participation of renewable producers reduces the clearing price on the demand curve. A key observation from Fig. \ref{fig:pjm_bid} is that the right half of the demand curve is well approximated by a linear function. We assume that the demand is at $1$ p.u. before renewable resources are introduced in the system. Therefore the clearing price decreases as a linear function of amount of renewables in the system up to modest penetration levels. Let $\sum_{i=1}^N w_i$ be the total amount of renewable power that is bid into the day-ahead market.
 We parametrize the day-ahead price as a function of the renewable bids as 
\begin{equation} \label{eqn:pw}
p(w_1,\dots,w_N)=p\left(\sum_{i=1}^N w_i\right)=1-\alpha \sum_{i=1}^N w_i,
\end{equation}
where $\alpha$ is a parameter that controls how fast the price decreases. Note the parameter $\alpha$ in \eqref{eqn:pw} is typically larger than unity as a small amount of renewables could reduce the day-ahead price significantly.  For example, if there is 14\% penetration of wind in Fig. \ref{fig:pjm_bid}, then the clearing price would be reduced to 50\% of the original price, which gives $\alpha=3.5$. Figure \ref{fig:pjm_bid} is based on PJM data, but the shape of the curve, especially the sharp rise at the right of the curve is common among most electricity markets \cite{Weron06}. 

If the total bid $\sum_{i=1}^N w_i$ is larger than $\frac{1}{\alpha}$, then the approximation in \eqref{eqn:pw} breaks down. One possible way to resolve this issue is to threshold $p$ to $0$ when the total bid amount exceeds $\frac{1}{\alpha}$. However, since it is never in the interest of the producers to have a total bid of more than $\frac{1}{\alpha}$, negative prices would never arise in the day-ahead market. A more fundamental issue is that at higher levels of penetration, the nonlinearity in the curve in Fig \ref{fig:pjm_bid} becomes important to the behavior of the renewable generators in the day-ahead market. Therefore, our results hold under low to moderate penetration of renewables, and the nonlinearity needs to be accounted to extend the results to higher penetrations. 

\textbf{Remark:} A relevant question is whether conventional generator will change their bids in the presence of renewables. From PJM reports\cite{PJM13}, it seems that conventional generators already bid their true cost, and therefore would not change their bid (bidding lower make no economic sense and bidding higher decrease the chance that they are cleared).


\section{Coalition of Renewable Producers With Fixed Penalty} \label{sec:penalty}
Suppose that the renewable producers are divided into $K$ groups and let $\mc S_1,\dots,\mc S_K$ denote the groups. Two extreme examples serve as benchmarks throughout the paper: the grand coalition, where all producers are in one group; and individual competition, where the producers compete as single players. For a group $S_k$, its profit is defined as
\begin{subequations} \label{eqn:pi_S}
\begin{align}
\pi_{\mc S_k} &=p(w_1,\dots,w_N) \sum_{i \in \mc S_k} w_i - q \E\left[\left(\sum_{i \in \mc S_k} w_i-\sum_{i \in \mc S_k} W_i\right)^+\right] \\
& \stackrel{(a)}{=} (1-\alpha \sum_{l=1}^N w_l) \sum_{i \in \mc S_k} w_i - q\E\left[\left(\sum_{i \in \mc S_k} w_i-\sum_{i \in \mc S_k} W_i\right)^+\right],
\end{align}
\end{subequations}
where $(a)$ follows from the linear price model in \eqref{eqn:pw}. Comparing \eqref{eqn:pi} and \eqref{eqn:pi_S}, we see that essentially the bid $w_i$ of a single producer is replaced by the bid $\sum_{i \in \mc S_k} w_i$ of a group of producers. Since the penalty term in the expectation is not linear, the profit of a group is not simply the sum of its individual parts and this leads to the benefit of aggregation. To study the effect of aggregation in depth, we adopt the following stochastic model. 

\subsection{Stochastic Model}
Let $W$ be a positive random variable with mean $\mu$. This random variable can be interpreted as already incorporating the forecasting information, and $\mu$ can be thought as the \emph{day-ahead forecast}. Therefore, the distribution of $W$ is the conditional distribution of the forecast error conditioned on the day-ahead forecast.  If there are $N$ producers, let $W_1,\dots,W_N$ be drawn \emph{identically} from the distribution of $W$ and the output of the $i$'th producer is 
$ W_i/N$
and the total output of $N$ producers is
\begin{equation*}
\frac1N \sum_{i=1}^N W_i.
\end{equation*}
This stochastic model allows us to keep the mean of the total amount of renewable power in the system constant and focus on the effect of a large number of producers. The next two sections study the cases when the producers are independent and when they are correlated. We do not specify a particular distribution for $W_i$ since the results hold for a wide case of distributions. Note that since we assume that the random variables are identically distributed, they should be produced by the same type of source. For a mixture of renewable sources, e.g. wind and solar, the identical assumption may not hold. Algorithm 1 in the Appendix partially addresses the mixed case by outlining a procedure for selecting groups based only on their covariance matrices. 

Since each coalition can only offer a bid in quantity, the profit model in \eqref{eqn:pi_S} sets up a Cournot competition among the different groups \cite{Varian06}. In \cite{Zhang14}, we consider in depth the game theoretic questions such as the existence and properties of Nash equilibria. In particular, we show that under broad conditions, it is always in the benefit of producers to form groups.   In the current paper, we are less concerned with such questions, and simply note that the ISO can impose rules on the size of coalitions and we are interested in determining the \emph{optimal} size of the groups. 

The following definition states when we consider a set of groups to be optimal. 
\begin{defn}
Given a set of groups $\mc S_1,\dots,\mc S_K$ that form a partition of $\{1,\dots,N\}$, let $\{w_1,\dots,w_K\}$ be the set of bids at a Nash equilibrium under \eqref{eqn:pi_S}. We say the set of groups is optimal if $\sum_{k=1}^K w_k=\frac{1}{\alpha}$. 
\end{defn}
Under the price model, $1/\alpha$ is the maximum amount (measured as a fraction of the total demand) of renewables that can be injected into the system. At this bid value, the day-ahead price becomes 
\begin{equation*}
1-\alpha \sum_i w_i=1-\alpha \frac{1}{\alpha}=0.
\end{equation*} 
Any additional injection would cause the price to go negative, which means that our price model in \eqref{eqn:pw} breaks down.

%
%
%

\subsection{Independent Forecast Errors} \label{sec:iid}
In this subsection, we will focus on the situation where the forecast errors of the renewable producers are independent. Of course, the forecast errors are correlated in the practice. However, understanding the independent case helps to illustrate the core concepts associated with grouping of producers. The main result is that neither individual competition nor the grand coalition are efficient, but there exists coalitions of intermediate size that are efficient.  The correlated forecast error case is considered in the next subsection using NREL data. 

\begin{figure}[!t]
\hspace{-0.4cm}
\subfigure[N=50]{
\includegraphics[scale=0.35]{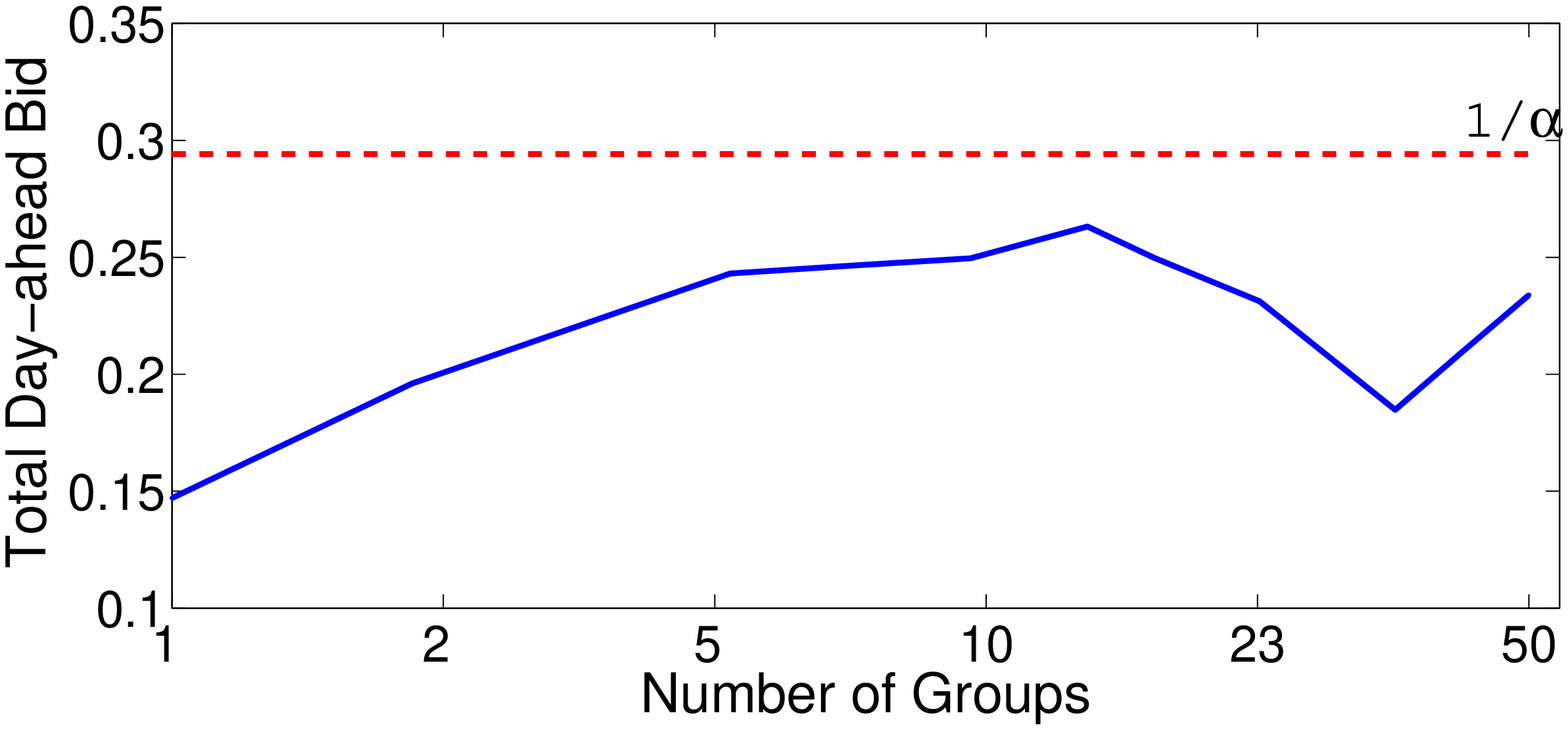}}
\subfigure[N=100]{
\hspace{-0.4cm}
\includegraphics[scale=0.35]{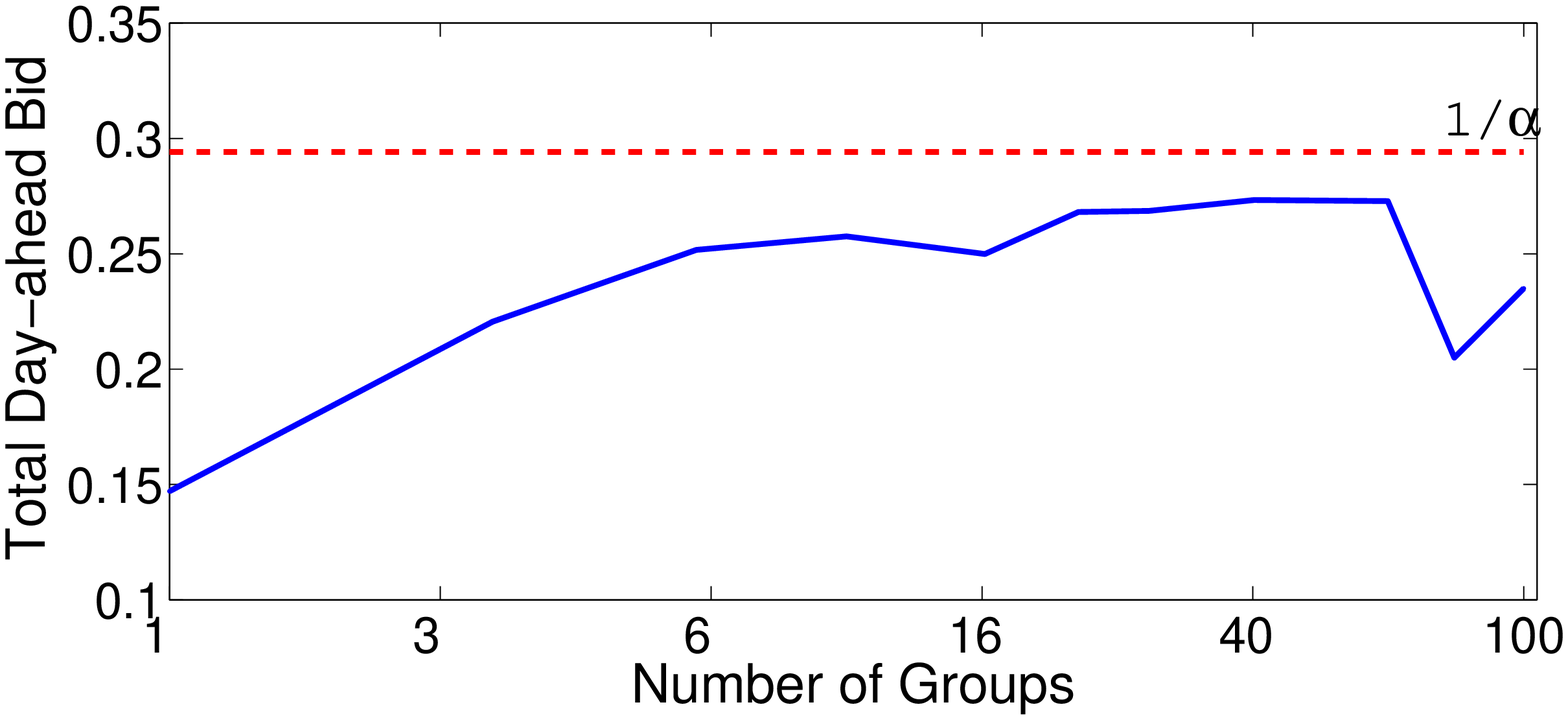}}
\subfigure[N=500]{
\hspace{-0.4cm}
\includegraphics[scale=0.38]{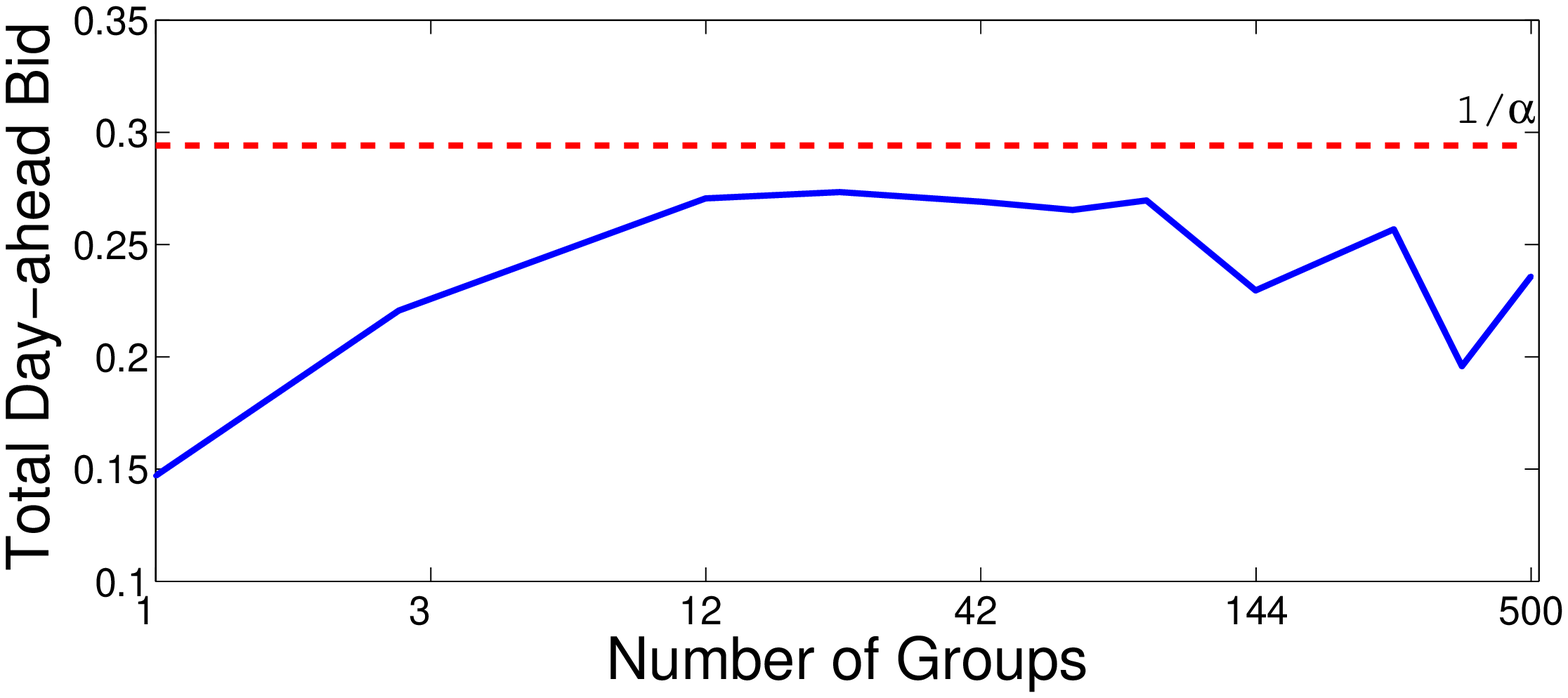}}
\subfigure[N=1000]{
\hspace{-0.4cm}
\includegraphics[scale=0.35]{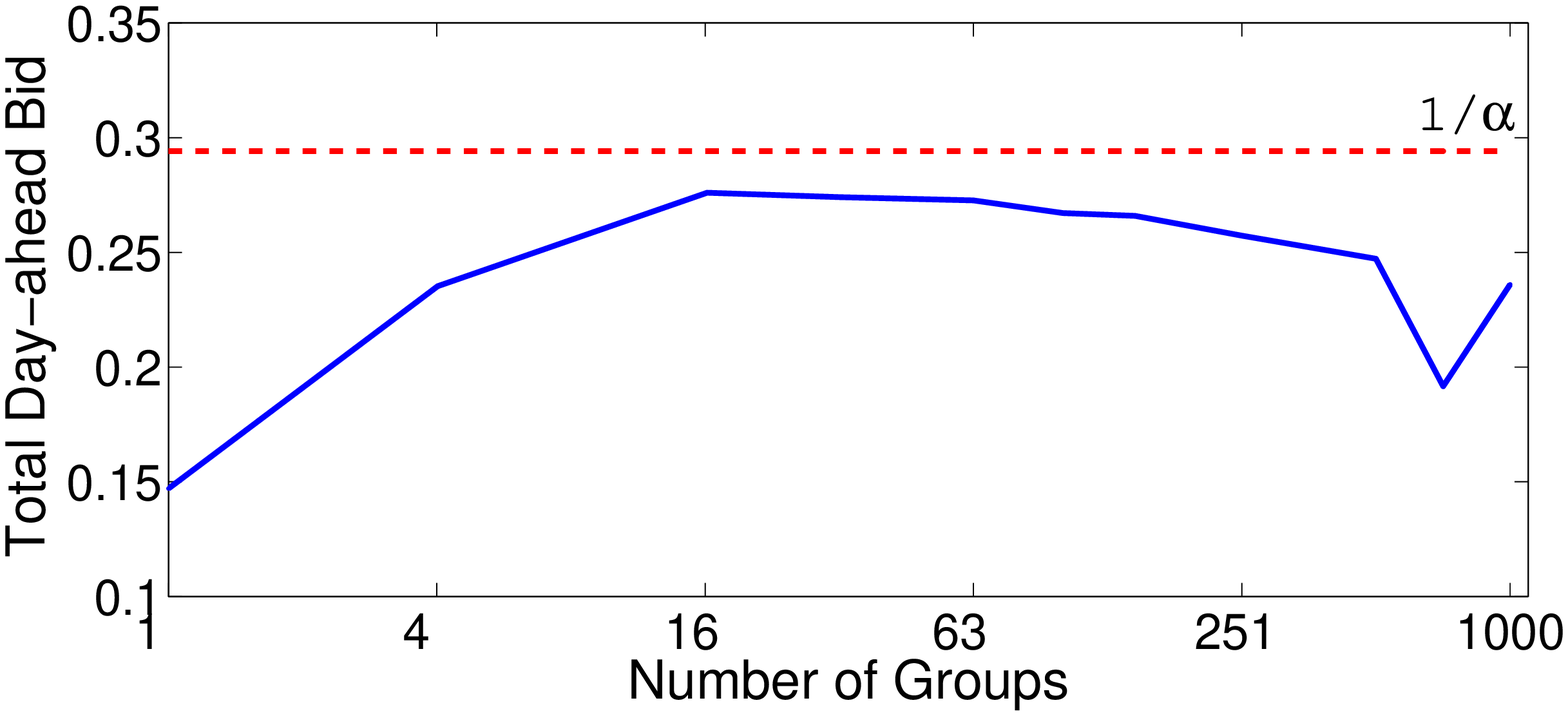}}
\caption{The total day-ahead bid versus the number of groups (in log-scale) for four different values of $N$. The leftmost point of the each plot is the grand coalition and the rightmost point is individual competition. We see that the grand coalition induces the least amount of biding in the system under the simulation parameters of $\alpha=0.3$ and $q=1$. The maximum amount of day-ahead bid always occurs at intermediate number of groups. As $N$ increase from 50 to 1000, the maximum total bid approaches $1/\alpha$, the efficient outcome. }
\label{fig:N}
\end{figure}

We consider coalitions of three types: the grand coalition, individual competition and groups of intermediate size. Of these three, we demonstrate that the groups of intermediate size is optimal from the system point of view since they balance the trade-off between market power and uncertainty. Figure \ref{fig:N} plots the total day-ahead bid versus the number of groups for four increasing values of $N$.
The leftmost point in the figures represent the grand coalition (a single group) and the rightmost point represent individual competition ($N$ groups). As shown in Fig. \ref{fig:N}, the maximum total bid occurs at group sizes in between the two extremes. As $N$ grows, the maximum approaches the $1/\alpha$ limit. 

The parameters used in generating Fig. \ref{fig:N} are $\alpha=3.4$, $q=1$ and $E[W_i]=\mu=0.3$. The value of $\mu$ specifies that the total expected renewable in the system is $30\%$. The $W_i$'s are assumed to be Gaussian (Gaussian forecast errors).  The exact value of $\mu$, $\alpha$, $q$ and the particular distribution of $W_i$'s do not change the qualitative behavior of the results. Below we explain why the intermediate case is optimal while the two extremes are suboptimal.


First we consider the case of the grand coalition. The profit of the coalition is given by
\begin{equation}
\pi(w)= (1-\alpha w) w- q \E\left[\left(w-\sum_{i=1}^N W_i/N\right)^+\right],
\end{equation}
where $w$ is the bid of the grand coalition. 
The bid $w^*$ that maximize the profit solves 
\begin{equation*}
1-2 \alpha w^* -q \Pr\left(w^* - \sum_{i=1}^N W_i/N\right)=0,
\end{equation*}
where $\Pr()$ is with respect to the joint probability of the $W_i$'s. We are interested in $w^*$ as $N$ increases. As $N$ grows, by the law of large numbers, $\sum_{i=1}^N W_i/N$ approaches its mean in probability. In this regime, there is essentially no uncertainty for the grand coalition and market power dominates. Figure \ref{fig:N} shows that as expected, the amount of renewable injected into the system is limited by this market power.  

Next, we consider the case where the individual producers do not form groups. The profit for producer $i$ is (special case of \eqref{eqn:pi_S})
\begin{equation} \label{eqn:pi_i} 
\pi_i (w_i)= \left(1- \alpha \sum_{l=1}^N w_l\right) w_i - q\E[(w_i-W_i/N)^+].
\end{equation}
Given the bid of other producers, the bid that maximizes the profit of producer $i$ is the solution to
\begin{equation}
1-2 \alpha w_i^* - \alpha \sum_{l \neq i}^N w_l - q \Pr(W_i/N \leq w_i^*)=0.
\end{equation}
Due to symmetry among the players, all producers would submit the same bid at equilibrium\footnote{Since the producers are symmetric, it is easy to show that if not all bids are the same, then some producers would change their bids, so the Nash equilibrium of this game is symmetric and unique.}. Therefore the optimal bid $w_i^*$ solves
\begin{equation*}
1-(N+1) \alpha w_i^*-q\Pr(W_i/N \leq w_i^*)=0.
\end{equation*}
As $N$ increases, each producer becomes smaller in size so they act as price takers. However, the second stage penalty dominates, so the total amount of renewable injection is still limited as shown in Fig. \ref{fig:N}. 

Groups of intermediate size balance the trade-off between market power and uncertainty.  As $N$ grows, the total number of groups also grows, ensuring competitiveness. On the other hand, the number of producers in a group also increases, ensuring averaging to reduce uncertainty. As Fig. \ref{fig:N} shows, this intermediate grouping maximizes the amount of renewable injected into the system.

\textbf{Remark: Optimal Group Size?} It turns out that most groups are asymptotically efficient, in the sense as long as the number of groups ($K$) and the number of producers in a group ($N/K$) both grows as $N$ grows. However, although the number of producers in practice is large, it is not infinite\footnote{Large control regions such as CAISO or ERCOT could have thousands of renewable producers.}. Therefore an interesting question is to find the \emph{optimal scaling rate}, or the group size that approaches $1/\alpha$ the fastest. As we show in  \cite{Zhang14}, the optimal scaling rate is obtained by $N^{2/3}$ groups of size $N^{1/3}$.  This result is obtained by separating the effect of market power and uncertainty in analyzing the efficiency of the Cournot game. 

\textbf{Remark: Producer Revenue} An interesting question about coalition formation is how individual producer's revenue change as the group size grows. For simplicity, consider the case with no uncertainty, that is, there is no second stage penalty. Let there be $N$ producers and they are divided into $K$ groups. It is straightforward to calculate that the per producer profit is given by 
\begin{equation*}
\frac{1}{\alpha N} \frac{K}{(K+1)^2}.
\end{equation*}
The $1/(\alpha N)$ factor is a common term independent of $K$. Therefore as $K$ (the number of groups) increases, the per producer profit decreases roughly as $1/K$. On one hand, this is not surprising, since at an efficient Nash equilibrium, the marginal profit of a producer equal to its marginal production cost. Since wind power producers have zero marginal cost, each producer have zero profit at the efficient equilibrium. 

From the producers' point of view,  they have a great amount incentive to form the largest possible coalition. Under current regulations, producers are not allowed to form any kind of groups. Then if operators allowed coalition formations, but limited the maximum size of the group to be the size that induces the most renewable penetration, the producers would naturally form the desired groups. However, 
instead of regulations, we should seek a market mechanism for to limit the size of groups. This question is not considered in this paper but is an interesting direction for further studies. 

\subsection{Correlated Forecast Errors}
\begin{figure}[t!]
\centering
\subfigure[Forecast Wind Power]{
\includegraphics[scale=0.27]{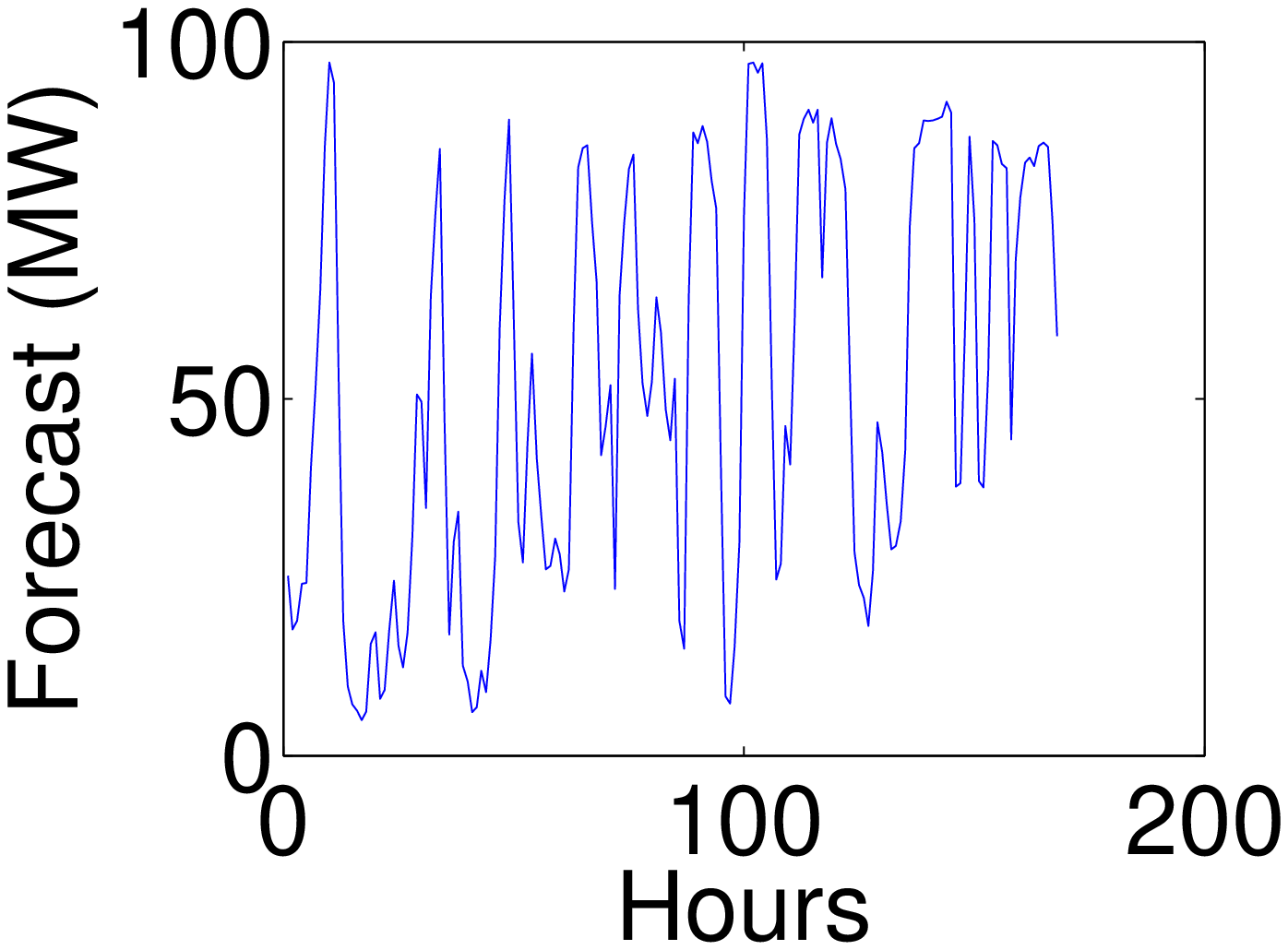}} 
\subfigure[Forecast Error]{
\includegraphics[scale=0.27]{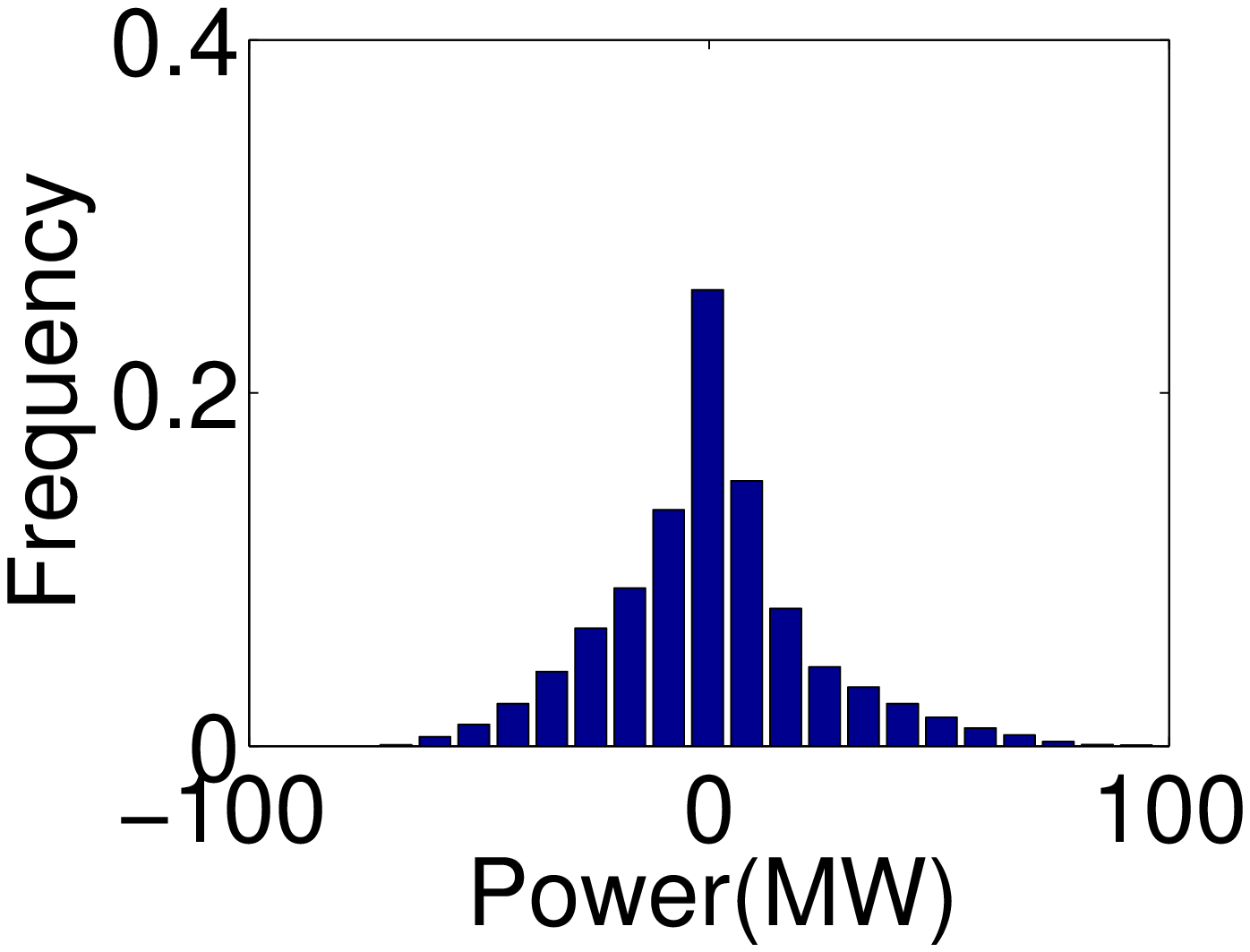}}
\caption{Forecasted wind power and its associated forecast error.}
\label{fig:error}
\end{figure}

In this section, we show that the coalitions of intermediate size are still optimal when the forecast errors are correlated. A general theory for correlated producers is difficult to develop in part due to the fact that results would depend on the particular distribution of the errors. In this section, we focus on empirical data from NREL eastern wind studies \cite{NRELeastern}. This dataset is a simulated study of the amount of wind power available at different geographical locations in the eastern part of U.S. Simulation was performed based on a meteorological and geographical conditions, validated using some field measurements. We consider the 302 locations that are in the PJM control area.  Figure \ref{fig:error} shows the forecasted wind power at a particular wind site and its associated forecast error. 

\begin{figure}[ht]
\centering
\includegraphics[scale=0.5]{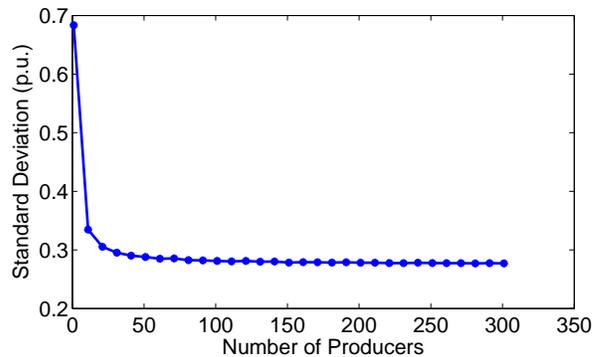}
\caption{Standard deviation as a function of the number of producers in an aggregation. The vertical axis is normalized by the total capacity of the aggregation.}
\label{fig:correlation}
\end{figure}

Figure \ref{fig:correlation} shows the standard deviation of the aggregate forecast ($1/N \sum_{i=1}^N W_i$) as a function of the number of wind farms in the aggregate. If the forecast errors are independent, we would expect that the standard deviation to decrease as the number of producers grows. Instead, Fig. \ref{fig:correlation} shows that the standard deviation flattens out as $N$ increases.

Similar to the previous section, we still look for a group structure that maximizes the total amount of wind power bid into the system. Note that the profit function for a group of producers is still given by \eqref{eqn:pi_S}. The main differences between the independent and correlated cases are: 1) the $1/\alpha$ result is not achievable since error does not approach zero as $N$ grows; and 2) groups size do need to not tend to infinity since all the benefit of averaging is achieved at finite group sizes.

For simplicity, we always normalize the wind farms to have equal capacity. However, since the standard deviations and the cross correlation between different producers are not equal, the optimal groups may not be of the same size. In fact, to find the best group structure is a combinatorial problem that requires the knowledge of the detailed joint distribution of all forecast errors. Rather than focus on finding the best possible way to group the wind farms,  we are more interested in the qualitative statement that some intermediate grouping is optimal for maximizing wind power injection. Therefore we present a simple greedy algorithm for selecting the groups. This algorithm  takes the number of groups as an input, and outputs a partition of the wind farms into the desired number of groups based only on the covariance matrix. It is given in the Appendix.
\begin{figure}[!t]
\hspace{-0.1cm}
\includegraphics[scale=0.5]{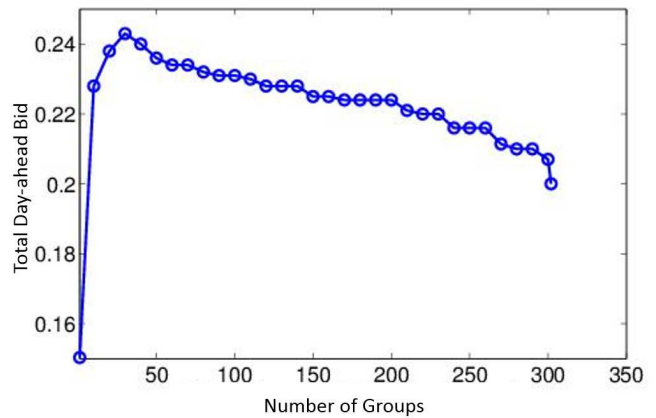}
\caption{Total amount of wind bid into the day-ahead market as a function of the number of groups for the wind farms in the NREL dataset. The maximum occurs at 30 groups, while the minimum occurs at a single group (grand coalition). Compare to Fig. \ref{fig:N}, the maximum is lower in the current case due to correlated errors.  }
\label{fig:correlated}
\end{figure}

Figure \ref{fig:correlated} plots the number of groups versus the total day-ahead bid of the groups. All the parameters are the same with respect to the simulations shown in Fig. \ref{fig:N}, in particular, $\alpha=3.4$. We see that in the correlated case, the maximum amount of wind power bid in day-ahead is about $24.3\%$ of the total demand, occurring when the producers are divided into 30 groups. Note that the total bid is less than the maximum amount in the i.i.d. (identically and independently distributed) case in the previous section, which is $1/\alpha=29\%$. This shows that not all the uncertainty can be averaged out of the system. Each of the groups contains roughly the same number of wind farms, ranging from 8 to 12 farms. The aggregate capacity of each group is also approximately the same. The approximate symmetry of the groups are expected since the capacity of each of the wind farms are normalized to be the same, and the original wind farms in the NREL dataset are similar to each other both in distribution and in capacity. The case of renewable producers that differ drastically in size is a direction of future research. 

\textbf{Remark: Market Power in Real-Time.} In Fig. \ref{fig:correlated}, the group selection is optimized by grouping producers with uncorrelated or anti-correlated errors to achieve uncertainty reduction. Because producers only exercise market power in the day-ahead stage, for groups of similar mean output, smaller uncertainty is preferred by both the system and the producers. However, if groups can exercise significant market power in real-time, positively correlated producers may group together to take advantage of large swings in their output. Exploring the tension created by real-time market power maybe a worthwhile question to answer. 

%

\section{Risk and Real-Time Prices} \label{sec:real-time}
In the model we have adopted so far, the real-time deviation is penalized when the bid quantity is less than the actual realized wind (recall \eqref{eqn:pi_i} and \eqref{eqn:pi_S}). In this section we show that the conclusion from previous sections is still valid under more complicated real-time mechanisms. In particular, we show that as long as the producers are risk averse, there is a benefit to forming groups of intermediate sizes. More broadly, as long as the objective functions for the producers are convex in $w_i-W_i$, a similar result can be derived as in Figs. \ref{fig:N} and \ref{fig:correlated}. Therefore the penalty adopted in \eqref{eqn:pi_i} and \eqref{eqn:pi_S} can be seen as a special case of a convex objective function. 

In practice, different market operators adopt different rules for handling real-time deviations. In some markets, both positive and negative deviations of generator output are penalized (e.g. Spain), although this double sided penalty is not applicable to renewable producers since their production can be easily curtailed. A more relevant real-time mechanism to consider is the real-time market. In many of the US electricity markets and some European markets (e.g. Nord Pool), a market is run at (or near to) real-time to readjust and balance supply and demand \cite{Pinson14}. A generator is then charged a penalty for its shortfall and charged or paid at the real-time market clearing price depending on the sign of its deviation.  Consider a group of renewable producers in $\mc S$. The profit of the group would be
\begin{equation} \label{eqn:rt_price} 
\begin{split}
\pi_S=& \left(1-\alpha \sum_{l=1}^N w_l \right) \sum_{i \in \mc S} w_i - \E \left[p_{\mbox{rt}} \cdot \left(\sum_{i \in \mc S} w_i-\sum_{i \in \mc S} W_i\right)\right] \\
& - q\E\left[\left(\sum_{i \in \mc S} w_i-\sum_{i \in \mc S} W_i \right)^+ \right]
\end{split}.
\end{equation}
The second term in \eqref{eqn:rt_price} represents the result of the settlement of the real-time market. If $p_{\mbox{rt}}$ is positive, and  $\sum_{i \in \mc S} w_i-\sum_{i \in \mc S} W_i$ is negative (excess power), the group receives payment by selling back its excess power. 

The modeling of the real-time price is notoriously difficult, in part due to its complex dependence on many possible aspects of the system (e.g. generator outages, unscheduled intertie flows, topology changes and others). For this section, we adopt the following model for the real-time price
\begin{equation} \label{eqn:rt}
\begin{split} 
p_{\mbox{rt}} &= \overline{p} \cdot \bd 1\left(\epsilon + \sum_{i}^N (w_i- W_i)\right) \\ &=
\begin{cases}
\overline{p} &\mbox{ if } \epsilon + \sum_{i}^N (w_i- W_i) >0 \\
0 & \mbox{ otherwise}
\end{cases}.
\end{split}
\end{equation}
In \eqref{eqn:rt}, the real-time price can take on two values, $\overline{p}$ and $0$. We think of $\overline{p}$ as the price cap of the system. The random variable $\epsilon$ is independent to the $W_i$'s and can be thought as an idiosyncratic shock in the demand. Together, the term $\epsilon + \sum_{i}^N (w_i- W_i)$ is interpreted as the net demand in the system at real-time. The real-time price hits the cap if the net demand is positive, and is zero otherwise. Of course, the model in \eqref{eqn:rt} is crude, but it does capture the volatility of the real-time prices \cite{Weron06}. 

To focus on the effect of real-time price on the profit of each group, we replace the penalty term in \eqref{eqn:rt_price} and consider a profit function of the form 
\begin{equation} \label{eqn:rt_only} 
\pi_S= \left(1-\alpha \sum_{l=1}^N w_l \right) \sum_{i \in \mc S} w_i - \E\left[p_{\mbox{rt}} \cdot \left(\sum_{i \in \mc S} w_i-\sum_{i \in \mc S} W_i \right)\right]. 
\end{equation}
As in Section \ref{sec:iid}, we focus on the case of i.i.d. $W_i$'s. Suppose there are $K$ equally sized groups $\mc S_1,\dots,\mc S_K$. The Nash equilibrium of bids is given by simultaneously finding  optimal solution to \eqref{eqn:rt_only} for each of the groups. Let $w_k$ be the equilibrium bid of $S_k$, and it is the solution of  
\begin{equation}\label{eqn:rt_K}
1-\alpha (K+1) w_k - \overline{p} \Pr\left(\sum_{k=1}^K- \epsilon < \sum_{k=1}^K w_k\right).
\end{equation}
Due to symmetry, the bids of all groups are equal. 

\begin{figure}[t!]
\centering
\includegraphics[scale=0.34]{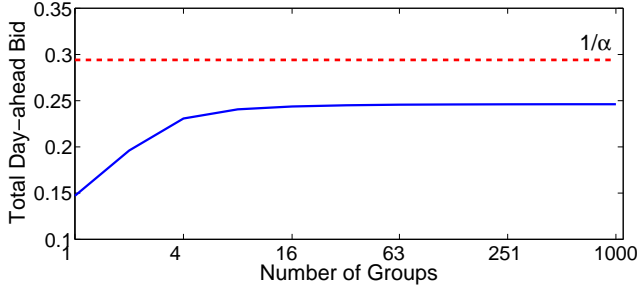}
\caption{The aggregate bid by the renewable producers as a function of the number of groups (log-scale) for 1000 producers. The grand coalition injects the least amount of renewable into the system and the bid increases as the number of groups increases. Note that it is not possible to approach  $1/\alpha$ because of the exogenous randomness in $\epsilon$.}
\label{fig:rt_N}s
\end{figure}

Figure \ref{fig:rt_N} shows the total day-ahead bid as a function of the number of groups for 1000 producers. We see that the grand coalition still bids the least amount of renewable power into the system and the total bid increases as the number of groups increases. In contrast to Fig. \ref{fig:N}, the maximum bid is achieved by individual competition.  Also, due to the presence of the exogenous random variable $\epsilon$, it is not possible to achieve an aggregate bid of $1/\alpha$. 

From Fig. \ref{fig:rt_N} it seems that grouping is not necessary since individual competition already induces maximum bidding among the producers. Indeed, this is not unexpected since the probability in \eqref{eqn:rt_K} depends only on the sum of all $W_i$'s and the exogenous random variable $\epsilon$. Therefore grouping does not affect the uncertainty part of the profit, and individual competition is optimal to maximize competition among the producers. Therefore, the \emph{mean} profit of a group is only affected by its market power, but not by the uncertainty of producers in that group.  

\begin{figure}[!t]
\subfigure[Individual Competition]{
\includegraphics[width=9cm]{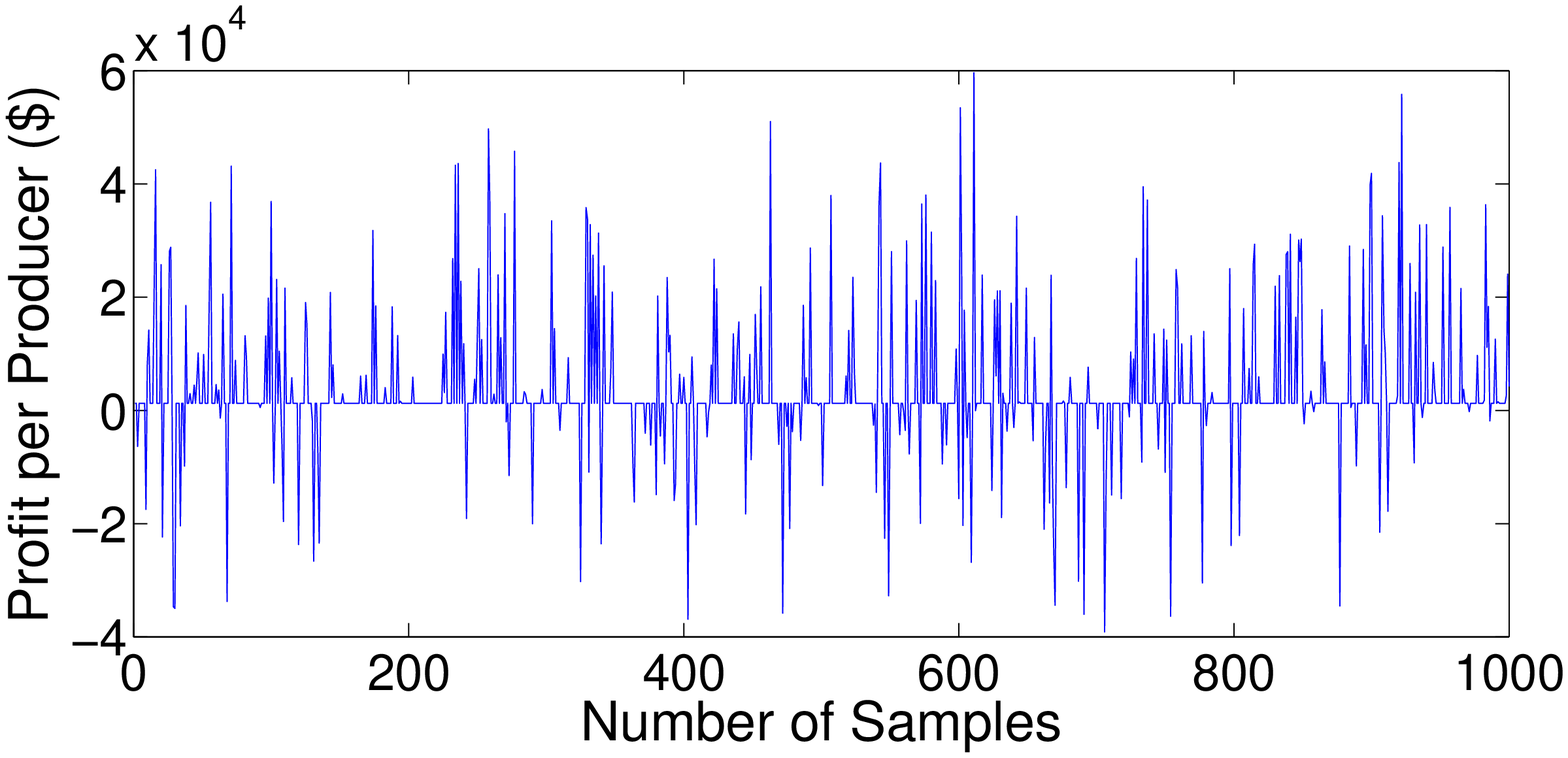}
\label{fig:rt_1}}
\subfigure[Groups of size 10]{
\includegraphics[width=9cm]{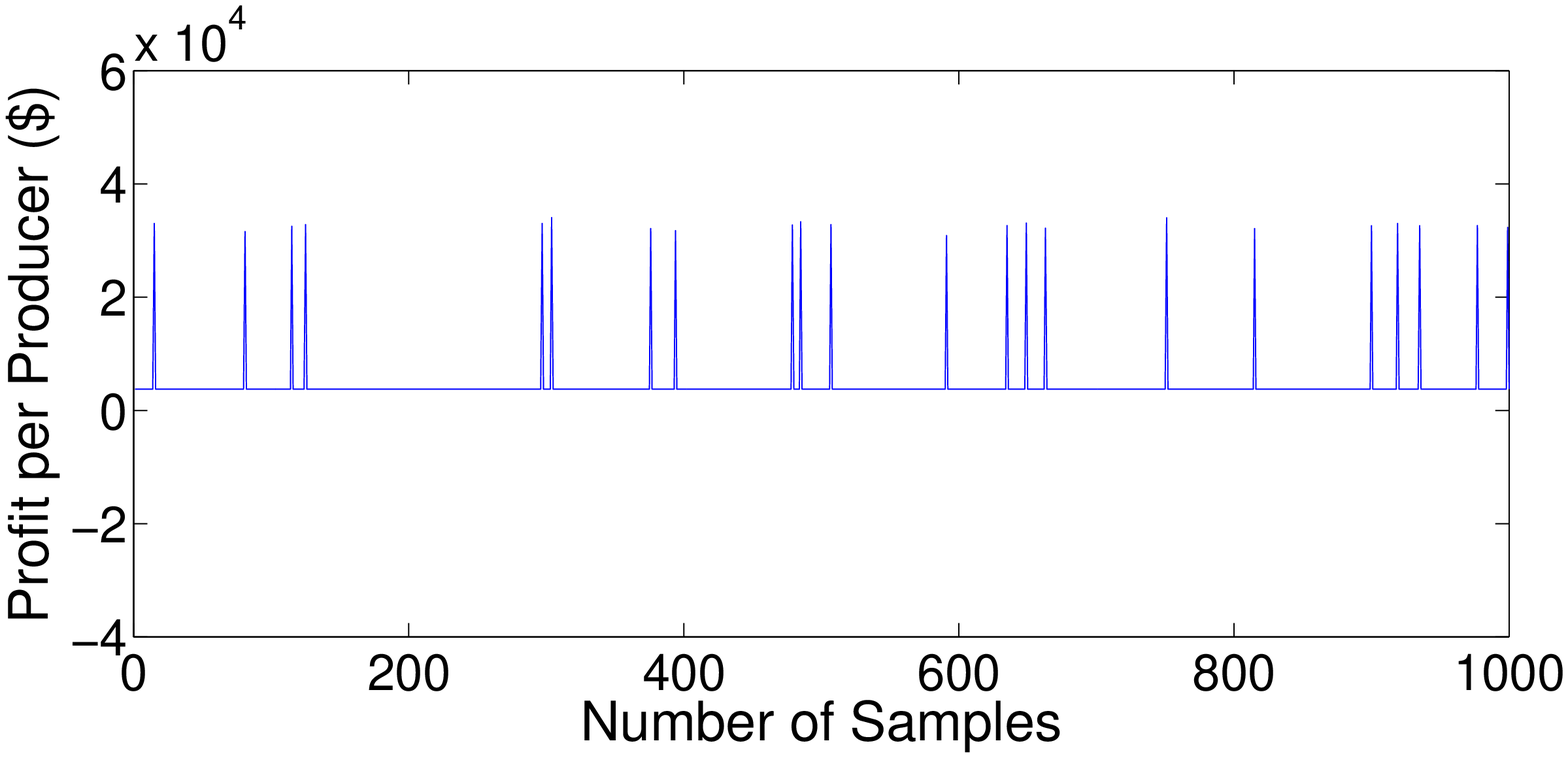}
\label{fig:rt_2}}
\caption{Sample paths of per producer profit for individual competition and groups of size 10. The total number of producer is 100. Two sample paths have the same mean, but the first is much more volatile than the second. }
\label{fig:rt_trace}
\end{figure}

Even though the average profit does not depend on uncertainty within a group, forming groups can still be beneficial to producers. As in Section \ref{sec:iid}, we focus on the case of i.i.d. $W_i$'s. Consider 100 producers ($N=100$). Assume $\epsilon$ is Gaussian with zero mean and standard deviation of $0.05$ (5\% uncertainty).  Figure \ref{fig:rt_trace} shows the sample paths of the per producer profit for individual competition and group of size 10. Both sample paths have the same mean value, meaning that the expected profit for each group is the same, but Fig.s \ref{fig:rt_1} is much more volatile than Fig. \ref{fig:rt_2}. From a producer's revenue management perspective, the revenue stream in Fig. \ref{fig:rt_2} is preferable. Therefore grouping of producers is still beneficial to the market, as long as the groups are not large enough to have significant market power. 

A form of risk not captured by the above model is the inexact knowledge of the probability distribution of $\epsilon_i$. Developing robust bidding solutions is an important direction for future work. 
  
\section{Conclusion}
In this paper we investigated group formations and strategic behaviors of renewable power producers in
electricity markets. To maximize the amount of renewable power injected into the system, we characterized the trade-off between market power and generation uncertainty as a function of the size of groups. We show there is a sweet spot in the sense that there exists groups that are large enough to achieve the uncertainty reduction of the grand coalition, but small enough such that they have no significant market power. We derived a linear model to quantify the effect of renewable bids on the day-ahead market clearing price. By modeling the day-ahead bidding process as a Cournot game, we  showed that grouping producers into coalitions of intermediate size achieve the social optimal outcome. We considered  both independent and correlated forecast errors under a fixed real-time penalty, and independent errors  under a real-time market where both selling and buying of energy is allowed.  We validated our claims using PJM and NREL data.
\appendix[Greedy Algorithm for Group Selection]
Given $N$ random variables $W_1,\dots,W_N$ and a positive integer $K<N$, we wish to partition the random variables into $K$ sets, $\mc S_1,\dots, \mc S_K$ that solves the following optimization problem
\begin{subequations} \label{eqn:sets}
\begin{align}
\max_{\mc S_1,\dots,\mc S_K} \; & \sum_{k=1}^K w_k \label{eqn:total_K} \\
\mbox{s.t. } & 1-2 \alpha w_k - \alpha \sum_{l \neq k} w_l - q \Pr(\sum_{i \in \mc S_k} W_i \leq w_k)=0 \; \forall \; k \label{eqn:nash_K}
\end{align},
\end{subequations}
where \eqref{eqn:total_K} is the total day-ahead bid and \eqref{eqn:nash_K} is the condition for the Nash equilibrium among the $K$ groups. 

The optimization in \eqref{eqn:sets} is a difficult problem for two reasons: it is combinatorial, and solving it requires the joint distribution of the random variables. To deal with the latter, we use the variance of  the aggregate $\sum_{i \in \mc S_k} W_i$ as a proxy for the uncertainty in $k$'th group; for the former, we adopt a greedy algorithm to construct the groups. 

The greedy algorithm is presented below. It proceeds in two parts. In the first part, it selects $K$ producers as seeds for the $K$ groups. In the second part, it progressively adds producers to the existing groups according to the covariance between them. In this paper, we compute empirical covariance between two wind farm's power output using their historical information, and use that in the algorithm. 

\begin{algorithm}
\renewcommand{\algorithmicrequire}{\textbf{Input:}}
\renewcommand{\algorithmicensure}{\textbf{Output:}}
\caption{Greedy Group Construction}\label{euclid}
\begin{algorithmic}
\REQUIRE $N >0$, $K < N$, random variables $W_1,\dots,W_N$
\ENSURE Partition of $1,\dots,N$ into sets $\mc S_1,\dots, S_K$
\STATE Let $\mc S_1,\dots,\mc S_K$ be empty sets 
\STATE Let $\mc S=\mc S_1 \cup \cdots \cup \mc S_K$ 
\STATE $\mc S_1 \leftarrow 1$
\STATE Define $\mc N ={2,\dots,N}$
\FOR{ $k=2$ \TO $K$ }
\STATE Find $i =   \arg\min_{i \in \mc N} \cov(W_i,\sum_{l \in \mc S})$
\STATE $\mc S_k \leftarrow i$
\STATE Remove $i$ from $\mc N$ 
\ENDFOR

\WHILE{$\mc N$ is not empty}
\STATE Let $i$ be an element in $\mc N$
\STATE Find $k^*=\arg\min_{k=1,\dots,K} \cov (W_i, \sum_{l \in \mc S_k} W_l)$
\STATE $\mc{S}_{k^*} \leftarrow i$
\STATE Remove $i$ from $\mc N$
\ENDWHILE
\end{algorithmic}
\end{algorithm}

\bibliography{mybib}
\bibliographystyle{IEEEtran}
\begin{IEEEbiography}[{\includegraphics[width=1in,height=1.25in,clip,keepaspectratio]{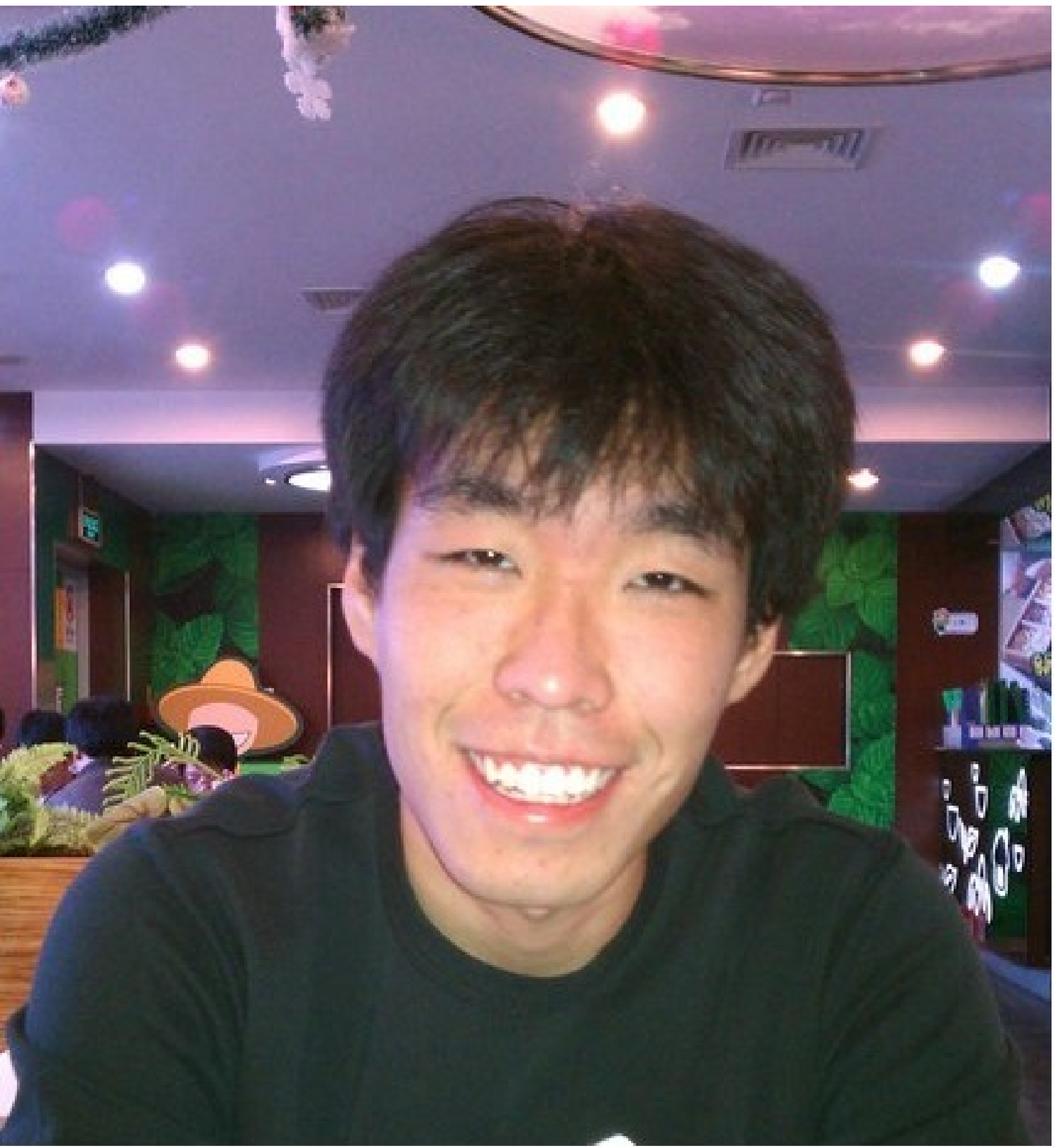}}]{Baosen Zhang} is a postdoctoral scholar at Stanford University, jointly hosted by departments of Civil and Environmental Engineering and Management \& Science Engineering. He will start as an Assistant Professor in Electrical Engineering at the University of Washington in 2015. 

He received the B.A.Sc. degree in engineering science from the University of Toronto, Toronto, ON, Canada, in 2008  and the Ph.D. degree Department of Electrical Engineering and Computer Science, University of California at Berkeley in 2013. His interest is in the area of power systems, particularly in the  fundamentals of power flow, data inference in SmartGrid and the economical challenges resulting from renewables. He was awarded several fellowships: the Post Graduate Scholarship from NSERC in 2011; the Canadian Graduate Scholarship from NSERC in 2008; the EECS fellowship from Berkeley in 2008.  
\end{IEEEbiography}

\begin{IEEEbiography}[{\includegraphics[width=1in,height=1.25in,clip,keepaspectratio]{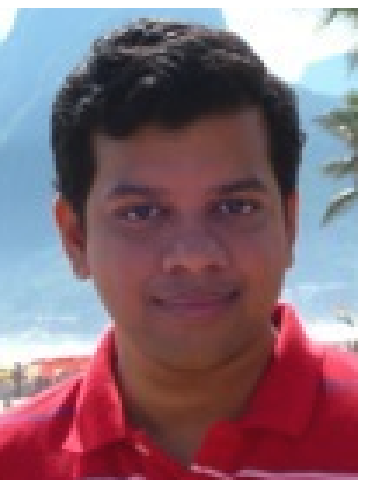}}]{Ram Rajagopal} received the B.Sc. degree in electrical
engineering from the Federal University of Rio
de Janeiro, Brazil; the M.Sc. degree in electrical and
computer engineering from the University of Texas,
Austin, TX, USA; and the M.A. degree in statistics
and the Ph.D. in electrical engineering and computer
sciences from the University of California, Berkeley,
CA, USA.

He is an Assistant Professor of Civil and Environmental
Engineering at Stanford University, Stanford,
CA, USA, where he directs the Stanford Sustainable
Systems Lab (S3L), focused on large scale monitoring, data analytics and stochastic
control for infrastructure networks.
His current research interests in power systems are in integration of renewables,
smart distribution systems and demand-side data analytics.
\end{IEEEbiography}

\begin{IEEEbiography}[{\includegraphics[width=1in,height=1.25in,clip,keepaspectratio]{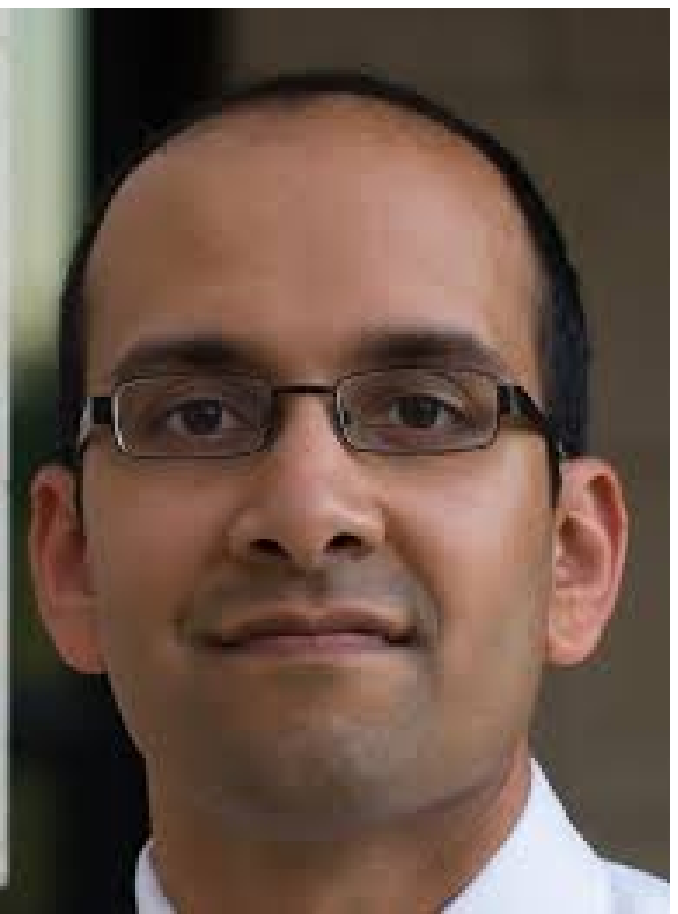}}]{Ramesh Johari} received the A.B. degree
in Mathematics from Harvard University (1998),
the Certificate of Advanced Study in Mathematics
from University of Cambridge (1999), and the Ph.D.
in Electrical Engineering and Computer Science
from M.I.T. (2004). He is currently an Associate
Professor of Management Science and Engineering,
and by courtesy, Computer Science and Electrical
Engineering, at Stanford University.
\end{IEEEbiography}

\end{document}